# Graph Irregularity Characterization with Particular Regard to Bidegreed Graphs

Ali Ghalavand,[1,*] Tamás Réti,[2] Igor Z. Milovanović[3] and Ali Reza Ashrafi[1]

[1]Department of Pure Mathematics, Faculty of Mathematical Sciences, University of Kashan, Kashan 87317-53153, I. R. Iran

[2]Donát Bánki Faculty of Mechanical and Safety Engineering, Óbuda University, Népszínház u. 8, H-1081, Budapest, Hungary

[3]Faculty of Electronic Engineering, Nis, Serbia

## Abstract

In this study we are interested mainly in investigating the relations between two graph irregularity measures which are widely used for structural irregularity characterization of connected graphs. Our study is focused on the comparison and evaluation of the discriminatory ability of irregularity measures called degree deviation $S(G)$ and degree variance $Var(G)$. We establish various upper bounds for irregularity measures $S(G)$ and $Var(G)$. It is shown that the Nikiforov's inequality which is valid for connected graphs can be sharpened in the form of $Var(G) < S(G)/2$. Among others it is verified that if G is a bidegreed graph then the discrimination ability of $S(G)$ and $Var(G)$ is considered to be completely equivalent.

## 1. Preliminary

Let $G$ be a connected simple graph for which $V(G)$, $E(G)$ denote the set of vertices and edges, and n and m the numbers of vertices and edges, respectively.

Let $\Delta = \Delta(G)$ and $\delta = \delta(G)$ be the maximum and the minimum degrees, respectively, of vertices of $G$. A universal vertex of an $n$-vertex graph $G$ is a vertex adjacent to all other vertices i.e. $\Delta(G) = n - 1$ holds. In a connected graph $G$ the number of universal vertices is denoted by $q = q(G)$. For a vertex $u$ in $G$, $d_G(u)$ stands for the degree of u. An edge of $G$ connecting vertices $u$ and $v$ is denoted by uv. It is assumed that the vertex sequence of G denoted by $(u_1, u_2, \ldots, u_n)$ is characterized by the degree sequence $\Delta = d_1 \geq d_2 \geq \ldots \geq d_n = \delta$ where $d_j$ denotes the degree of the vertex $u_j$. The degree set of a graph G is denoted by Ds(G).

---

[*] Corresponding author (Email: alighalavand@grad.kashanu.ac.ir).

For a connected graph G with n vertices and m edges the **average** degree is equal to 2m/n. The so-called cyclomatic number of a connected graph with n vertices and m edges is defined as $c = c(G) = m - n + 1$. A connected graph $G$ having $c(G) = k \geq 1$ cycles is said to be a k-cyclic graph. If $k = 0$ holds the corresponding acyclic graph is called a tree graph.

Using the standard terminology [1,2,3,4], a graph G is called R-regular if all its vertices have the same degree $R$. A connected graph is called irregular if it contains at least two vertices with different degrees. The number of vertices of degree $i$ is denoted by $N_i$.

An irregular graph $G$ is said to be $k$-degreed, if its degree set $Ds(G)$ contains exactly $k$ different degrees. Consequently, a *bidegreed graph* $G(\Delta, \delta)$ is an irregular graph whose vertices have exactly two different degrees.

An $n \geq 4$ vertex bidegreed graph $G(\Delta, \delta)$ is called a *balanced graph* if $n$ is an even integer for which $N_\Delta = N_\delta = n/2$ holds. Among $n$-vertex connected bidegreed graphs, path $P_4$ with degree sequence (2,2,1,1) and the so-called diamond graph $G_D$ with degree sequence (3,3,2,2) represent the smallest balanced bidegreed graphs.

## 2. Graph Irregularity Measures

By definition, a topological invariant $IT(G)$ is called an irregularity measure of a graph $G$ if $IT(G) \geq 0$ and $IT(G) = 0$ if and only if $G$ is a regular graph. The majority of irregularity measures are the so-called degree-based graph invariants.

The degree deviation $S(G)$ and the degree variance $Var(G)$ of a graph $G$ belong to the family of most popular graph irregularity measures. The degree deviation $S(G)$ was introduced by Nikiforov [5], and for a connected irregular graph $G$ with $n$ vertices and m edges it is defined by $S(G) = \sum_{i=1}^{n} \left| d_i - \frac{2m}{n} \right|$. The degree variance Var(G) introduced by Bell [6] is formulated as

$$Var(G) = \frac{1}{n} \sum_{i=1}^{n} \left( d_i - \frac{2m}{n} \right)^2 = \frac{M_1(G)}{n} - \left( \frac{2m}{n} \right)^2$$

where $M_1(G)$ stands for the first Zagreb index of a graph $G$. By definition [7, 8, 9, 10] $M_1(G) = \sum_{i=1}^{n} d_i^2$. Additionally, in our study we introduce two novel graph irregularity measures denoted by $IRD(G)$ and $IRR(G)$, both of them will be used in the subsequent sections. They are formulated as $IRD(G) = \frac{2N_\Delta N_\delta}{N_\delta + N_\Delta}(\Delta - \delta)$ and $IRR(G) = \frac{n}{2}(\Delta - \delta)$.

In this study we investigate primarily the relations among the above irregularity measures and evaluate their discriminating ability for characterizing the topological structure of connected bidegreed graphs. It should be mentioned that during last three decades the irregularity measures

$S(G)$ and $Var(G)$ and the corresponding extremal graphs have been intensively studied, some important results are included in [6, 11, 12, 13, 14].

## 3. Some Relationships between Irregularity Measures $S(G)$, $Var(G)$, $IRD(G)$ and $IRR(G)$

The aim of this section is to provide some inequalities between the irregularity measures $S(G)$, $Var(G)$, $IRD(G)$ and $IRR(G)$. In all inequalities, a necessary and sufficient condition for equality will be presented.

**Theorem 1.** [15] Let $x = (x_1, x_2, \ldots, x_n)$ be a real number sequence with the property $\sum_{i=1}^{n} x_i = 0$ and $\sum_{i=1}^{n} |x_i| = 1$. Then $|\sum_{i=1}^{n} a_i x_i| \leq \frac{1}{2}\left(\max_{1 \leq i \leq n} a_i - \min_{1 \leq i \leq n} a_i\right)$, where $a_1, a_2, \ldots, a_n$ are real numbers.

**Theorem 2.** If $G$ is an $n$-vertex graph with degree sequence $d(G) = (d_1, d_2, \ldots, d_n)$, then

(i) $Var(G) = \frac{1}{n}\sum_{i=1}^{n}\left(d_i - \frac{2m}{n}\right)^2 \leq \frac{\Delta - \delta}{2n} S(G)$,

(ii) $S(G) \leq \frac{n}{2}(\Delta - \delta) = IRR(G)$,

(iii) $Var(G) \leq \frac{1}{n^2}(IRR(G))^2 = \frac{(\Delta - \delta)^2}{4}$,

with equalities if and only if G is regular or a bidegreed balanced graph.

**Proof.** (i) Suppose $x_i = \frac{d_i - \frac{2m}{n}}{S(G)}$ and $a_i = \left|d_i - \frac{2m}{n}\right|$ for all $i = 1, 2, \ldots, n$. It is easy to see that $\sum_{i=1}^{n} x_i = 0$ and $\sum_{i=1}^{n} |x_i| = 1$. Apply Theorem 1 to deduce that

$$\left|\sum_{i=1}^{n} \frac{\left(d_i - \frac{2m}{n}\right)^2}{S(G)}\right| \leq \frac{1}{2}\max\left\{\Delta - \frac{2m}{n}, \frac{2m}{n} - \delta\right\} \leq \frac{1}{2}(\Delta - \delta).$$

Therefore,

$$\sum_{i=1}^{n}\left(d_i - \frac{2m}{n}\right)^2 \leq \frac{1}{2}(\Delta - \delta)S(G) \tag{1}$$

with equality if and only if G is regular or a bidegreed balanced graph.

(ii) By Cauchy-Schwarz inequality we have $\sum_{i=1}^{n}\left(d_i - \frac{2m}{n}\right)^2 \geq \frac{\left(\sum_{i=1}^{n}\left|d_i - \frac{2m}{n}\right|\right)^2}{\sum_{i=1}^{n} 1}$. This means that

$$\sum_{i=1}^{n}\left(d_i - \frac{2m}{n}\right)^2 \geq \frac{S(G)^2}{n} \tag{2}$$

with equality if and only if $\left|d_1 - \frac{2m}{n}\right| = \left|d_2 - \frac{2m}{n}\right| = \cdots = \left|d_n - \frac{2m}{n}\right|$ if and only if $|Ds(G)| \leq 2$ and $n\left(\Delta - \frac{2m}{n}\right) = S(G) = n\left(\frac{2m}{n} - \delta\right)$ if and only if $|Ds(G)| \leq 2$, $N_\delta = n - N_\Delta$ and $n(\Delta + \delta) = 4m$ if and only if $|Ds(G)| \leq 2$ and $(n - 2N_\delta)(\Delta - \delta) = 0$ if and only if $G$ is regular or a bidegreed balanced graph. Thus, by (1) and (2), $S(G) \leq \frac{n}{2}(\Delta - \delta) = IRR(G)$, with equality if and only if G is regular or a bidegreed balanced graph.

(iii) By (ii) and (1), $Var(G) \leq \frac{1}{n^2}(IRR(G))^2$ with equality if and only if G is regular or a bidegreed balanced graph. $\square$

**Proposition 3.** [16]: For connected bidegreed graphs $G(\Delta, \delta)$ with $n$ vertices and $m$ edges it has been proved that $S(G(\Delta, \delta)) = IRD(G(\Delta, \delta))$.

Suppose G is a graph and k is a positive integer. Define $N(G, k) = \{v \in V(G) | d_G(v) = k\}$.

**Theorem 4.** Let $G$ be graph with $n$ vertices. Then $S(G) \geq \frac{2N_\Delta N_\delta}{n}(\Delta - \delta) = \frac{N_\Delta + N_\delta}{n} IRD(G)$ and the equality holds if and only if $G$ is regular or a bidegreed graph.

**Proof.** By definition of $N(G, k)$,

$$S(G) = \sum_{v \in V(G)} \left|d_G(v) - \frac{2m}{n}\right| \geq \sum_{v \in N(G,\delta)} \left(\frac{2m}{n} - \delta\right) + \sum_{v \in N(G,\Delta)} \left(\Delta - \frac{2m}{n}\right)$$

$$= N_\delta \left(\frac{2m}{n} - \delta\right) + N_\Delta \left(\Delta - \frac{2m}{n}\right). \tag{3}$$

The equality holds in (3) if and only if $|Ds(G)| = 1$, $|Ds(G)| = 2$ or $(|Ds(G)| = 3$ and $\frac{2m}{n} \in D(G))$. On the other hand,

$$\frac{2m}{n} - \delta = \frac{2m - n\delta}{n}$$

$$\geq \frac{N_\Delta(\Delta - \delta)}{n}, \tag{4}$$

$$\Delta - \frac{2m}{n} = \frac{n\Delta - 2m}{n} \geq \frac{N_\delta(\Delta - \delta)}{n}. \tag{5}$$

The equality in (4) (and also in (5)) holds if and only if $|Ds(G)| \leq 2$. Now, By Equations (3–5),

$S(G) \geq \frac{N_\delta N_\Delta(\Delta-\delta)}{n} + \frac{N_\Delta N_\delta(\Delta-\delta)}{n} \geq \frac{2N_\Delta N_\delta}{n}(\Delta - \delta) = \frac{N_\Delta + N_\delta}{n} IRD(G)$ and the equality holds if and only if $|Ds(G)| \leq 2$. $\square$

**Theorem 5.** If G is an n-vertex graph, then $Var(G) \geq \frac{(N_\delta + N_\Delta)^2}{n^4} IRR(G) IRD(G)$. The equality is satisfied if and only if $|Ds(G)| \leq 2$.

**Proof.** By definition of $N(G, k)$,

$$Var(G) = \frac{1}{n}\sum_{v\in V(G)}\left(d_G(v) - \frac{2m}{n}\right)^2 \geq \frac{1}{n}\sum_{v\in N(G,\delta)}\left(\frac{2m}{n} - \delta\right)^2 + \frac{1}{n}\sum_{v\in N(G,\Delta)}\left(\Delta - \frac{2m}{n}\right)^2$$

$$= \frac{1}{n}N_\delta\left(\frac{2m}{n} - \delta\right)^2 + \frac{1}{n}N_\Delta\left(\Delta - \frac{2m}{n}\right)^2$$

with equality if and only if $|Ds(G)| \leq 2$. So, by (4) and (5),

$$Var(G) \geq \frac{N_\delta N_\Delta(N_\delta + N_\Delta)(\Delta - \delta)^2}{n^3} = \frac{(N_\delta + N_\Delta)(\Delta - \delta)}{2n^4}IRD(G)$$

$$= \frac{(N_\delta + N_\Delta)^2}{n^4}IRR(G)IRD(G)$$

with equality if and only if G is $|Ds(G)| \leq 2$.  ▫

**Corollary 6.** Let $G$ be an $n$-vertex connected bidegreed graph. Then

1. $Var(G) = \frac{N_\delta + N_\Delta}{n^2}(\Delta - \delta)^2$ for any bidegreed graph.

2. $IRD(G) = S(G) \leq IRR(G)$ with equality if and only if G is a bidegreed balanced graph.

3. $Var(G) = \frac{1}{n^2}IRR(G)IRD(G) \leq \frac{1}{n^2}(IRR(G))^2 = \frac{(\Delta-\delta)^2}{4}$ with equality if and only $G$ is a bidegreed balanced graph.

**Proof.** Since G is bidegreed, $n = N_\delta + N_\Delta$ and the proof follows from Theorems 2, 4 and 5.  ▫

**Proposition 7.** ([17]) Let $G(\Delta, \delta)$ be a connected bidegreed graph with $n\ (\geq 4)$-vertices. Then $Var(G(\Delta, \delta)) \leq \frac{1}{4}(\Delta - \delta)^2$ with equality if and only if $G(\Delta, \delta)$ is a balanced bidegreed graph.

Proof. According to Corollary 6, $Var(G(\Delta, \delta)) = \frac{N_\Delta N_\delta}{n^2}(\Delta - \delta) = \frac{N_\Delta(n-N_\Delta)}{n^2}(\Delta - \delta)^2$. Consider the monotonically increasing function $h(N_\Delta)$ defined in the interval $\left(2 \leq N_\Delta \leq \frac{n}{2}\right)$, $h(N_\Delta) = \frac{N_\Delta(n-N_\Delta)}{n^2} = \frac{N_\Delta}{n} - \frac{N_\Delta^2}{n^2}$. From $\frac{\partial h(N_\Delta)}{\partial N_\Delta} = \frac{1}{n} - \frac{2}{n^2}N_\Delta \geq 0$ it follows that $h(N\Delta)$ has a maximum value when $N_\Delta = n/2$ holds. Consequently, $Var(G(\Delta, \delta)) \leq \frac{1}{4}(\Delta - \delta)^2$ where equality holds if and only if $G(\Delta, \delta)$ is a balanced bidegreed graph.  ▫

**Proposition 8.** Let $G = G(\Delta, \delta)$ be a connected bidegreed graph with $n$-vertices and m-edges. Then $S(G(\Delta, \delta)) = IRD(G(\Delta, \delta)) \leq IRR(G(\Delta, \delta))$, where equality holds if and only if $G(\Delta, \delta)$ is a bidegreed balanced graph (i.e. if $n \geq 4$ is an even integer, and $N_\Delta = N_\delta = \frac{n}{2}$).

**Proof.** By Proposition 3, $S(G(\Delta, \delta)) = IRD(G(\Delta, \delta)) = \frac{2N_\Delta(n-N_\Delta)}{n}(\Delta - \delta)$. Consider the monotonically increasing function $g(N_\Delta)$ defined in the interval $(2, n/2)$ by $g(N_\Delta) = \frac{2N_\Delta N_\delta}{n} = \frac{2N_\Delta(n-N_\Delta)}{n}$. Its maximum value with respect to $N_\Delta$ can be computed from $\frac{\partial g(N_\Delta)}{\partial N_\Delta} = \frac{2}{n}(n - 2N_\Delta) \geq 0$. As can be observed, function $g(N_\Delta)$ has a maximum value if $N_\Delta = N_\delta = n/2$ is fulfilled.

Consequently, for $n$-vertex bidegreed graphs $S(G(\Delta, \delta)) = IRD(G(\Delta, \delta)) \leq IRR(G(\Delta, \delta))$, with equality if and only if $n \geq 4$ is an even integer and $N_\Delta = N_\delta = n/2$. □

**Remark 1.** If $n$ is an odd integer, then for any $n$-vertex bidegreed graph $S(G) < IRR(G)$ holds. Moreover, if $n$ is even integer, but $N_\Delta \neq N_\delta$ then for a bidgreed graph, $IRR(G)$ will be always larger than $S(G)$. For example, if $G$ is the 4-vertex star $K_{1,3}$ then $S(K_{1,3}) = 3 < 4 = IRR(K_{1,3})$.

## 4. Upper Bounds for $S(G)$ and $Var(G)$ Measures

A connected irregular graph $G$ having at least one universal vertex is called a *dominating graph*. Dominating graphs form a subset of irregular graphs of diameter 2. There are several bidegreed graphs belonging to the family of dominating graphs. Typical examples are the $n \geq 5$ vertex wheel graphs $W_n$ and the complete split graphs [3, 13, 14].

As it is known, an $n$-vertex complete split graph denoted by $CS(n, k)$ consists of an independent set of $n - k$ vertices and a clique of $k$ universal vertices, such that each vertex of the clique is adjacent to each vertex of the independent set. [13, 14]. Complete split graphs $CS(n, k)$ represent a particular subclass of bidegreed dominating graphs $G(\Delta, \delta)$ where $\Delta = n - 1$ and $\delta = k$, and the corresponding edge number is equal to $m(CS(n, k)) = [(2n - 1)k - k^2]/2 = [(2n - 1)\delta - \delta^2]/2$.

**Lemma 9.** [18] Let $G$ be an $n$-vertex connected irregular graph. If $M_1(G)$ is maximal among $n$-vertex irregular graphs then $G$ contains at least one universal vertex of degree $n - 1$.

**Lemma 10.** [19] Let $G$ be a connected graph with $n$ vertices and $m$ edges. Then $M_1(G) \leq \frac{2m[2m+(n-1)(\Delta-\delta)]}{n+\Delta-\delta}$ with equality if $G$ is a regular graph or $G$ is isomorphic to the complete split graph $CS(n, k)$ where $\delta = k$ and $m = [(2n - 1)\delta - \delta^2]/2$.

**Proposition 11.** Let $G$ be a dominating graph with $n$-vertices and $m$-edges. Then

$$Var(G) = \frac{M_1(G)}{n} - \left(\frac{2m}{n}\right)^2 \leq \left(\frac{2m}{n}\right)\frac{[2m + (n-1)(n-1-\delta)]}{2n - 1 - \delta} - \left(\frac{2m}{n}\right)^2$$

and equality holds if $G$ is isomorphic to a complete split graph $CS(n, k)$ where $\delta = k$ and $m = [(2n - 1)\delta - \delta^2]/2$.

Proof. Because for complete split graphs with $n$ vertices and m edges $\Delta = n - 1$, $\delta = k$, and $m = [(2n - 1)\delta - \delta 2]/2$ hold, the result follows from Lemma 10.

**Proposition 12.** [5, 17]: For a connected graph $G$ with $n$ vertices and $m$ edges, it is verified that

$$S(G) \leq n\sqrt{Var(G)} \tag{6}$$

In what follows, we show that there are infinitely many balanced bidegreed graphs for which equality holds in Eq. (6).

A broad class of bidegreed balanced graphs can be constructed from an n-vertex R≥3 regular graph $G_R$ by inserting vertices of degree 2 into $G_R$. The concept of construction is illustrated with graphs depicted in Fig.1.

**Example 1.** Balanced bidegreed graphs of different type can be generated from the 3-regular polyhedral graph $G_T$. It is the graph of the 6-vertex trigonal prism. Balanced bidegreed graphs $G_j$ ($j = 1, 2, 3$) constructed from $G_T$ are 12-vertex non-isomorphic graphs, all of them contain 6 vertices of degree 3 and 6 vertices of degree 2.

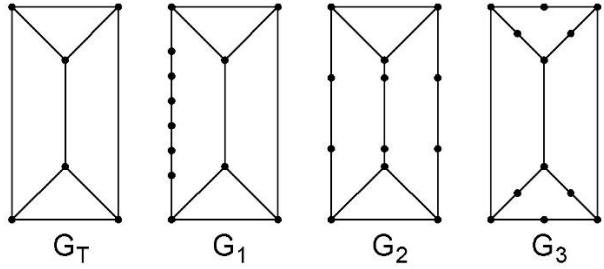

**Figure 1.** Balanced bidegreed graphs $G_j$ ($j = 1,2,3$) constructed from the 3-regular graph $G_T$.

Balanced bidegreed graphs $G_j$ ($j = 1, 2, 3$) have identical degree deviation $S(G_j) = IRD(G_j) = IRR(G_j) = 6$ and identical degree variation $Var(G_j) = \frac{1}{12}\sum_{i=1}^{12}\left(d_i - \frac{30}{12}\right)^2 = \frac{1}{4}$. As can be observed, balanced bidegreed graphs depicted in Fig. 1 cannot be discriminated by irregularity measures S(G) and Var(G).

**Proposition 13.** Let $G(\Delta, \delta)$ be an $n$-vertex balanced bidegreed graph. Then $S(G(\Delta, \delta)) = n\sqrt{Var(G(\Delta, \delta))} = \frac{n}{2}(\Delta - \delta) = IRR(G(\Delta, \delta))$.

Proof. For balanced bidegreed graph $G(\Delta, \delta)$ the equality $N_\Delta = N_\delta = n/2$ holds. It follows that $\frac{2m}{n} = \frac{\Delta N_\delta + \delta N_\delta}{n} = \frac{\Delta+\delta}{2}$. On one hand, $S(G) = N_\Delta\left(\Delta - \frac{\Delta+\delta}{2}\right) + N_\delta\left(\frac{\Delta+\delta}{2} - \delta\right) = \frac{n}{2}\left\{\frac{2\Delta-\Delta-\delta}{2} + \frac{\Delta+\delta-2\delta}{2}\right\} = \frac{n}{4}\{\Delta - \delta + \Delta - \delta\} = \frac{n}{2}(\Delta - \delta)$. Also, $n \times Var(G) = N_\Delta\left(\Delta - \frac{\Delta+\delta}{2}\right)^2 + N_\delta\left(\frac{\Delta+\delta}{2} - \delta\right)^2 = \frac{(\Delta-\delta)^2}{4}$. Consequently, $S(G(\Delta, \delta)) = n\sqrt{Var(G(\Delta, \delta))} = \frac{n}{2}(\Delta - \delta) = IRR(G(\Delta, \delta))$.

## 5. Comparing the Discrimination Power of Irregularity Measures S(G) and Var(G)

A common property of irregularity measures $Var(G)$ and $S(G)$ is that if two n-vertex graphs $G_A$ and $G_B$ have an identical degree sequence, then $Var(G_A) = Var(G_B)$ and $S(G_A) = S(G_B)$ hold.

During last three decades the discrimination power of various graph irregularity measure was evaluated in several relavant publications [20−26].

In what follows it will be demonstrated for bidegreed graphs that there is a strict algebraic relation between irregularity measures $S(G)$ and $Var(G)$. This observation implies that for an n-vertex bidegreed graph the discrimination performance of $Var(G)$ and $S(G)$ is considered to be completely equivalent.

**Proposition 14.** Let $G(\Delta, \delta)$ be an $n$-vertex connected bidegreed graph. Then $Var\big(G(\Delta, \delta)\big) = \frac{\Delta-\delta}{2n} S(G(\Delta, \delta))$.

**Proof.** According to Corollary 6, for the bidegreed graph $G(\Delta, \delta)$, $S\big(G(\Delta, \delta)\big) = IRD(G) = \frac{2N_\Delta N_\delta}{n}(\Delta - \delta)$ and $Var\big(G(\Delta, \delta)\big) = \frac{1}{n^2} IRR(G) IRD(G) = \frac{N_\Delta N_\delta}{n^2}(\Delta - \delta)^2$ hold. From these formulas the result follows.

Using comparative tests, it is revealed that if G is not bidegreed, then for general irregular graphs there is no uniquely defined deterministic algebraic relation between irregularity measures $S(G)$ and $Var(G)$.

In Fig. 2 this fact is demonstrated by irregular graphs with 6 vertices and 12 edges. For irregular graphs shown in Fig. 2 the corresponding topological indices are summarized in Table 1.

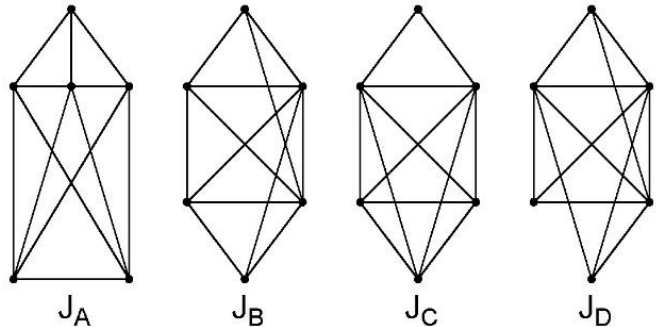

**Figure 2.** All connected irregular graphs with 6 vertices and 12 edges.

**Table 1.** Topological invariants of irregular graphs with 6-vertices and 12-edges

| Graph | Topological parameter | | | |
|---|---|---|---|---|
| | $M_1$ | S | Var | q |
| $J_A$ | 98 | 2 | 1/3 | 1 |
| $J_B$ | 100 | 4 | 2/3 | 2 |
| $J_C$ | 102 | 4 | 1 | 2 |
| $J_D$ | 102 | **6** | 1 | 3 |

As we can conclude, graphs $J_A$, $J_B$ and $J_C$ are tridegreed graphs, while $J_D$ is a balanced bidegreed graph having 3 vertices of degree 3 and 3 vertices of degree 5.

Among all irregular graphs depicted in Fig.2, there exist irregular graphs with q = 1,2 and 3 universal vertices. The bidegreed graph $J_D$ with 3 universal vertices is the only graph with a maximum irregularity measure S, i.e. $S(J_D) = 6$.

Table 1 contains the corresponding first Zagreb indices and the numbers q of universal vertices for graphs $J_A$, $J_B$, $J_C$ and $J_D$.

As we have already mentioned that Cioabă proved [19] that if among n-vertex connected graphs, the first Zagreb index $M_1(G)$ is maximal for a graph $G$, then $G$ contains at least one universal vertex of degree $n-1$. Consequently, if $G$ is an $n$-vertex connected graph for which $M_1(G)$ is maximal, then $G$ contains at least one universal vertex.

Let $G_{n,m}$ be a connected irregular graph with n vertices and m edges. From the definition of Var(G) it follows that the value of $Var(G_{n,m})$ will be maximal if $M_1(G_{n,m})$ is maximal, as well.

In Table 1 graphs $J_C$ and $J_D$ have equal maximal first Zagreb index, i.e. $M_1(J_C) = M_1(J_D) = 102$ holds. It follows from this observation that the maximal degree variance (Var=1) belongs to graphs $J_C$ and $J_D$. Moreover, we can conclude that the maximal degree deviation of S=6, belongs to bidegreed graph $J_D$.

In order to compare different irregularity measures IA(G) and IB(G) for connected irregular graphs we introduced a topological invariant ω(G) formulated as $\omega(G) = \frac{IA(G)}{IB(G)}$. As a particular case, consider the graph invariant defined by $\Omega(G) = \frac{Var(G)}{S(G)}$.

Nikiforov mentioned without proof in Ref. [5] that if $G$ is a connected irregular graph then inequality $\Omega(G) = \frac{Var(G)}{S(G)} \leq 1$ holds. In what follows, we find a better bond for $\Omega(G)$.

**Corollary 15.** $\Omega(G) < \frac{1}{2}$.

**Proof.** From inequality (1), one can see that $\Omega(G) = \frac{Var(G)}{S(G)} \leq \frac{1}{2n}(\Delta - \delta) \leq \frac{n-2}{2n} = \frac{1}{2} - \frac{1}{n} < \frac{1}{2}$. ☐

The inequality $0 \leq Var(G) \leq S(G)/2$ is the best possible, it cannot be sharpened as it is demonstrated by the following two examples.

**Example 2.** Consider the infinite sequence of wheel graphs $W_n$ with $n \geq 5$ vertices. For these bidegreed graphs $\Delta - \delta = n - 4$ holds. Because

$$Var(W_n) = \frac{(n-1)(n-4)^2}{n^2} \to \infty \text{ if } n \to \infty,$$

$$S(W_n) = \frac{2(n-1)(n-4)}{n} \to \infty \text{ if } n \to \infty.$$

One obtains that $\Omega(W_n) = \frac{Var(W_n)}{S(W_n)} = \frac{n-4}{2n} = \frac{1}{2} - \frac{2}{n} \to \frac{1}{2}$ if $n \to \infty$.

**Example 3.** There exists an infinite sequence of $n$-vertex paths $P_n$ for which $Var(P_n) = \frac{2(n-2)}{n^2} \to 0$ if $n \to \infty$, $S(P_n) = \frac{4(n-2)}{n} \to 4$ if $n \to \infty$, it follows that $\Omega(P_n) = \frac{Var(P_n)}{S(P_n)} \to 0$ if $n \to \infty$.

Starting with considerations outlined previously the following proposition yields.

**Proposition 16.** For a connected irregular graphs $\Omega(G) < 1/2$ and $Var(G) < S(G)/2$ are fulfilled. Moreover, for any small $\varepsilon > 0$ positive number there exists an $n$-vertex graph $G_n$ for which $\left|\frac{Var(G_n)}{S(G_n)} - \frac{1}{2}\right| < \varepsilon$ holds if n is sufficiently large. For example, such graphs are the $n$-vertex stars $K_{1,n-1}$ and $n$-vertex wheel graphs $W_n$.

**Proposition 17.** Let $H_1(\Delta_1, \delta_1)$ and $H_2(\Delta_2, \delta_2)$ be $n$-vertex connected bidegreed graphs for which $DIF = \Delta_1 - \delta_1 = \Delta_2 - \delta_2 \geq 1$ holds. Then

$$\Omega(H_1) = \frac{Var(H_1)}{S(H_1)} = \Omega(H_2) = \frac{Var(H_2)}{S(H_2)} = \frac{DIF}{2n}.$$

**Proof.** Because $H_1(\Delta_1, \delta_1)$ and $H_2(\Delta_2, \delta_2)$ are connected bidegreed graphs with $n$ vertices, it follows that for $= 1, 2$ $\Omega(H_i) = \frac{\Delta_j - \delta_j}{2n} = \frac{DIF}{2n}$. ▨

**Example 4.** Consider the connected bidegreed 7-vertex graphs depicted in Fig. 3 They have identical $DIF = \Delta - \delta$ parameter and different edge numbers ($m = 6, 8, 10$ and $16$).

For path P$_7$ the equality $\Delta - \delta = 2 - 1 = 1$ holds, and for bicyclic and 4-cyclic graphs B$_7$ and F$_7$ one obtains that $\Delta - \delta = 3 - 2 = 1$, respectively. Graph G$_7$ is a 10-cyclic bidegreed graph for which $\Delta - \delta = 5 - 4 = 1$. Because for these 7-vertex graphs $DIF = \Delta - \delta = 1$ holds, we have that $\Omega(P_7) = \Omega(B_7) = \Omega(F_7) = \Omega(G_7) = \frac{1}{14}$.

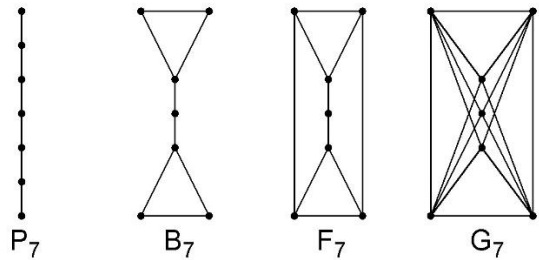

**Figure 3.** Bidegreed 7-vertex graphs with different edge numbers.

## 6. Multiplicative decomposition of degree variance

**Proposition 18.** [28] Let $G$ be a connected graph with $n$ vertices and $m$ edges. Then

$$Var(G) = \frac{1}{n}\sum_{i=1}^{n}\left(d_i - \frac{2m}{n}\right)^2 = \frac{M_1(G)}{n} - \left(\frac{2m}{n}\right)^2 \leq \left(\Delta - \frac{2m}{n}\right)\left(\frac{2m}{n} - \delta\right).$$

Based on the previous proposition, for connected bidegreed graphs a multiplicative decomposition of the degree variance can be performed.

**Proposition 19.** Let $G(\Delta, \delta)$ be a connected bidegreed graph with $n$ vertices and $m$ edges. Then

$$Var(G(\Delta, \delta)) = \frac{1}{n}\sum_{i=1}^{n}\left(d_i - \frac{2m}{n}\right)^2 = \frac{M_1(G(\Delta, \delta))}{n} - \left(\frac{2m}{n}\right)^2 = \left(\Delta - \frac{2m}{n}\right)\left(\frac{2m}{n} - \delta\right).$$

**Proof.** It follows that $\left(\Delta - \frac{2m}{n}\right)\left(\frac{2m}{n} - \delta\right) = \frac{2m}{n}(\Delta + \delta) - \Delta\delta - \left(\frac{2m}{n}\right)^2$. Consequently, it suffices to show that $M_1(G(\Delta, \delta)) = 2m(\Delta + \delta) - n\Delta\delta$. Because $n = N_\Delta + N_\delta$, $2m = \Delta N_\Delta + \delta N_\delta$ and $M_1(G(\Delta, \delta)) = \Delta^2 N_\Delta + \delta^2 N_\delta$, one obtains that $(\Delta N_\Delta + \delta N_\delta)(\Delta + \delta) - (N_\Delta + N_\delta)\Delta\delta = \Delta^2 N_\Delta + \delta^2 N_\delta = M_1(G(\Delta, \delta))$. ▣

It is interesting to note that a similar formula can be derived for the so-called 2-walk linear graphs representing a broad class of connected irregular graphs [29-34].

A connected irregular graph $G$ is called 2-walk linear (more precisely, 2-walk $(a, b)$– linear) if there exists a unique rational number pair (in fact, integers pair) $(a, b)$ such that $S(u) = d_G(u)t(u) = ad_G(u) + b$ holds for every vertex $u$ of $G$, where $S(u)$ is the sum of the degrees of all vertices adjacent to $u$, and $t(u)$ stands for the average of degrees of vertices adjacent to u.

As it is known, [29, 30, 31, 32, 33]:

i) If $G$ is a 2-walk $(a, b)$ linear graph, then $a$ and $b$ are integers, and $a \geq 0$.

ii) A connected irregular graph has exactly two main eigenvalues $\lambda > \mu$ if and only if it is 2-walk linear.

iii) If $G$ is a 2-walk $(a, b)$ linear graph with two main eigenvalues $\lambda$ and $\mu$, then $\lambda, \mu = \frac{1}{2}(a \pm \sqrt{a^2 + 4b})$, consequently $\lambda + \mu = a$ and $\lambda\mu = -b$ where $\lambda$ is the spectral radius of $G$.

**Proposition 20.** [34]: Let $G$ be a 2-walk $(a, b)$ linear graph with main eigenvalues $\lambda$ and $\mu$, where $\lambda > \mu$. Then $Var(G) = \left(\lambda - \frac{2m}{n}\right)\left(\frac{2m}{n} - \mu\right)$.

i) As an example, consider the 11-vertex and 20-edge connected graph GR with degree sequence $(5^1, 4^5, 3^5)$. It is depicted in Fig.4. This triangle-free graph having chromatic number 4 is known as the Grötzsch graph [4]. Graph $GR$ is a 2-walk $(a,b)$ linear graph with parameters $a = 1$ and $b = 10$.

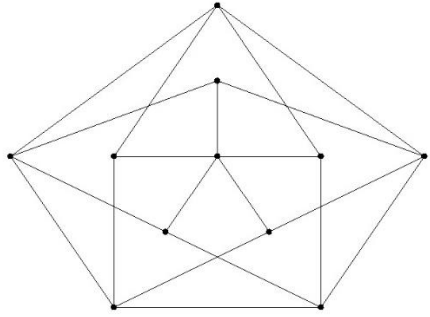

**Figure 4.** A 11-vertex 2-walk $(a,b)$ linear graph.

The corresponding main eigenvalues of graph $GR$ are $\lambda(GR), \mu(GR) = \frac{1}{2}(1 \pm \sqrt{41})$. It is easy to check that $Var(G) = \left(\lambda(GR) - \frac{40}{11}\right)\left(\frac{4}{11} - \mu(GR)\right) = \frac{50}{121}$.

**Remark 2.** There are several bidegreed 2-walk linear graphs with various structure, such graphs are the wheel graphs, friendship graphs, complete split graphs.

## 7. $S(G)$ and $Var(G)$ for n-vertex connected graphs with cyclomatic number ≤ (n + 2)/2

The aim of this section is to study the graph invariants $S(G)$ and $Var(G)$ for $n$-vertex connected graphs with cyclomatic number c such that $1 \leq c \leq \frac{n+2}{2}$.

**Theorem 21.** [35] If $G$ is a connected graph with $n$ vertices and cyclomatic number $c$, then

$$N_1 = 2 - 2c + \sum_{i=3}^{\Delta}(i-2)N_i \text{ and } N_2 = 2c + n - 2 - \sum_{i=3}^{\Delta}(i-1)N_i,$$

where $N_i$ for all $i = 1,2,\ldots,\Delta$ are the number of vertices of degree $i$ in $G$.

**Theorem 22.** Let $T$ be a tree on $n$ vertices. Then,

1. $S(T) = \frac{4(n-2)}{n} + \frac{2(n-2)}{n}\sum_{i=3}^{\Delta} N_i(i-2)$,

2. $Var(T) = \frac{2(n-2)}{n^2} + \frac{1}{n}\sum_{i=3}^{\Delta} N_i(i-1)(i-2)$,

3. $IRR(T) = \frac{n}{2}(\Delta - 1)$,

4. $IRD(T) \leq \frac{4N_\Delta(\Delta-1) + 2N_\Delta \sum_{i=3}^{\Delta}(i-2)N_i(\Delta-1)}{2+(\Delta-1)N_\Delta}$, with equality if and only if $D(G) \subseteq \{1,2,\Delta\}$.

**Proof.** Since $m = n - 1$, by definition of $S(T)$ and Theorem 21 we have

$$S(T) = N_1\left(\frac{2(n-1)}{n} - 1\right) + N_2\left(2 - \frac{2(n-1)}{n}\right) + \sum_{i=3}^{\Delta} N_i\left(i - \frac{2(n-1)}{n}\right)$$

$$= 2\left(\frac{2(n-1)}{n} - 1\right) + \sum_{i=3}^{\Delta}(i-2)N_i\left(\frac{2(n-1)}{n} - 1\right) + n\left(2 - \frac{2(n-1)}{n}\right)$$

$$- 2\left(2 - \frac{2(n-1)}{n}\right) - \sum_{i=3}^{\Delta}(i-1)N_i\left(2 - \frac{2(n-1)}{n}\right) + \sum_{i=3}^{\Delta} N_i\left(i - \frac{2(n-1)}{n}\right)$$

$$= \frac{4(n-2)}{n} + \frac{2}{n}\sum_{i=3}^{\Delta} N_i(i-2)(n-2).$$

This completes (1). To prove (2), we apply definition of $Var(T)$ and Theorem 21 to deduce that

$$n \times Var(T) = N_1\left(\frac{2(n-1)}{n} - 1\right)^2 + N_2\left(2 - \frac{2(n-1)}{n}\right)^2 + \sum_{i=3}^{\Delta} N_i\left(i - \frac{2(n-1)}{n}\right)^2$$

$$= 2\left(\frac{2(n-1)}{n} - 1\right)^2 + \sum_{i=3}^{\Delta}(i-2)N_i\left(\frac{2(n-1)}{n} - 1\right)^2 + n\left(2 - \frac{2(n-1)}{n}\right)^2$$

$$- 2\left(2 - \frac{2(n-1)}{n}\right)^2 - \sum_{i=3}^{\Delta}(i-1)N_i\left(2 - \frac{2(n-1)}{n}\right)^2$$

$$+ \sum_{i=3}^{\Delta} N_i\left(i - \frac{2(n-1)}{n}\right)^2 = \frac{2(n-2)}{n} + \sum_{i=3}^{\Delta} N_i(i-1)(i-2).$$

Hence, $Var(T) = \frac{2(n-2)}{n^2} + \frac{1}{n}\sum_{i=3}^{\Delta} n_i(i-1)(i-2)$ which completes (2). The proof of (3) is straightforward and so it is omitted. To prove (4), we use definition of $IRD(T)$ and Theorem 21. We can see that,

$$IRD(T) = \frac{2N_1 N_\Delta(\Delta - 1)}{N_1 + N_\Delta} = \frac{4N_\Delta(\Delta - 1) + 2N_\Delta \sum_{i=3}^{\Delta}(i-2)N_i(\Delta - 1)}{2 + \sum_{i=3}^{\Delta}(i-2)N_i + N_\Delta}$$

$$\leq \frac{4N_\Delta(\Delta - 1) + 2N_\Delta \sum_{i=3}^{\Delta}(i-2)N_i(\Delta - 1)}{2 + (\Delta - 2)N_\Delta + N_\Delta}$$

$$= \frac{4N_\Delta(\Delta - 1) + 2N_\Delta \sum_{i=3}^{\Delta}(i-2)N_i(\Delta - 1)}{2 + (\Delta - 1)N_\Delta},$$

with equality if and only if $Ds(G) \subseteq \{1,2,\Delta\}$. ⬜

**Corollary 23.** If $T$ is a tree on $n$ vertices, then

1. $S(T) \geq \frac{4(n-2)}{n}$ and equality holds if and only if $T \cong P_n$.

2. $Var(T) \geq \frac{2(n-2)}{n^2}$ and equality holds if and only if $T \cong P_n$.

3. $IRR(T) \geq \frac{n}{2}$ and equality holds if and only if $T \cong P_n$.

4. $IRD(T) \geq \frac{4N_\Delta(\Delta-1)+2N_\Delta^2(\Delta-2)(\Delta-1)}{2+(\Delta-2)(n-N_1-N_2)+N_\Delta}$, with equality if and only if $Ds(G) \subseteq \{1,2,\Delta\}$.

5. $IRD(T) \leq \frac{4N_\Delta(\Delta-1)+2N_\Delta(\Delta-2)(\Delta-1)(n-N_1-N_2)}{2+(\Delta-1)N_\Delta}$, with equality if and only if $Ds(G) \subseteq \{1,2,\Delta\}$.

**Proof.** The proofs of Parts (1), (2), (3) and (5) are direct consequences of Theorem 22. To prove (4), we have

$$IRD(T) = \frac{4N_\Delta(\Delta-1)+2N_\Delta \sum_{i=3}^{\Delta}(i-2)N_i(\Delta-1)}{2+\sum_{i=3}^{\Delta}(i-2)N_i + N_\Delta}$$

$$\geq \frac{4N_\Delta(\Delta-1)+2N_\Delta \sum_{i=3}^{\Delta}(i-2)N_i(\Delta-1)}{2+(\Delta-2)(n-N_1-N_2)+N_\Delta}$$

$$\geq \frac{4N_\Delta(\Delta-1)+2N_\Delta^2(\Delta-2)(\Delta-1)}{2+(\Delta-2)(n-N_1-N_2)+N_\Delta},$$

with equality if and only if $Ds(G) \subseteq \{1,2,\Delta\}$. ▫

**Corollary 24.** If $T$ is a tree on $n$ vertices, then

$$Var(T) - \frac{1}{2n}S(T) = \frac{1}{n}\sum_{i=3}^{\Delta} N_i(i-2)\left(i - \frac{2n-2}{n}\right).$$

**Proof.** The proof follows from Theorem 22. ▫

**Corollary 25.** If $T$ is a tree on $n \geq 4$ vertices, then $\Omega(T) = \frac{Var(T)}{S(T)} \geq \frac{1}{2n}$ and equality holds if and only if $T \cong P_n$.

**Proof.** By Corollary 24, $Var(T) - \frac{1}{2n}S(T) \geq 0$, equality holds if and only if $T$ is isomorphic to $P_n$. Therefore, $\Omega(T) = \frac{Var(T)}{S(T)} \geq \frac{1}{2n}$ and equality holds if and only if $T \cong P_n$. ▫

**Corollary 26.** Let $T$ be a tree on $n$ vertices. Then $S(T) \geq IRD(G)$, equality holds if and only if $T$ is bidegreed tree.

**Proof.** If tree $T$ is bidegreed, then $T = T(\Delta \geq 2, \delta = 1)$, and according to Proposition 3 the equality $S(T) = IRD(T)$ holds. Assuming that $T$ is a non-bidegreed tree with $\Delta \geq 3$, then by Theorem 22 we have

$$S(T) - IRD(T)$$
$$\geq \frac{4(n-2)}{n} + \frac{2}{n}\sum_{i=3}^{\Delta} N_i(i-2)(n-2)$$
$$- \frac{4N_\Delta(\Delta-1) + 2N_\Delta \sum_{i=3}^{\Delta}(i-2)N_i(\Delta-1)}{2 + (\Delta-1)N_\Delta}$$
$$= \frac{8(n-(\Delta-1)N_\Delta - 2)}{n(2+(\Delta-1)N_\Delta)} + \sum_{i=3}^{\Delta} \frac{(4i-8)(n-(\Delta-1)N_\Delta - 2)}{n(2+(\Delta-1)N_\Delta)}$$
$$\geq \frac{8N_2}{n(n-N_2)} + \sum_{i=3}^{\Delta} \frac{(4i-8)N_2}{n(n-N_2)} \geq \frac{8N_2}{n(n-N_2)} + \frac{(4\Delta-8)N_2 N_\Delta}{n(n-N_2)} \geq 0,$$

with equality if and only if $Ds(T) = \{1, \Delta\}$. ▫

**Theorem 27.** Let $G$ be a connected graph with $n$ vertices and cyclomatic number $c$ such that $1 \leq c \leq \frac{n+2}{2}$. Then,

$$S(G) = \frac{2}{n}\sum_{i=3}^{\Delta} N_i (in - 2m),$$

$$Var(G) = \frac{1}{n}\sum_{i=3}^{\Delta} N_i(i-1)(i-2) - \frac{2(2m-n)(m-n)}{n^2}.$$

**Proof.** Since $m = n + c - 1$ and $1 \leq c \leq \frac{n+2}{2}$, $2 \leq \frac{2m}{n} \leq 3$. Therefore, by Theorem 21 we have

$$S(G) = N_1\left(\frac{2m}{n} - 1\right) + N_2\left(\frac{2m}{n} - 2\right) + \sum_{i=3}^{\Delta} N_i\left(i - \frac{2m}{n}\right)$$

$$= (2-2c)\left(\frac{2m}{n} - 1\right) + \sum_{i=3}^{\Delta}(i-2)N_i\left(\frac{2m}{n} - 1\right)$$

$$+ (2c+n-2)\left(\frac{2m}{n} - 2\right) - \sum_{i=3}^{\Delta}(i-1)N_i\left(\frac{2m}{n} - 2\right) + \sum_{i=3}^{\Delta} N_i\left(i - \frac{2m}{n}\right)$$

$$= \frac{2}{n}\sum_{i=3}^{\Delta} N_i (in - 2m),$$

$$n \times Var(G) = N_1 \left(\frac{2m}{n} - 1\right)^2 + N_2 \left(\frac{2m}{n} - 2\right)^2 + \sum_{i=3}^{\Delta} N_i \left(i - \frac{2m}{n}\right)^2$$

$$= (2 - 2c)\left(\frac{2m}{n} - 1\right)^2 + \sum_{i=3}^{\Delta}(i - 2)N_i \left(\frac{2m}{n} - 1\right)^2$$

$$+ (2c + n - 2)\left(\frac{2m}{n} - 2\right)^2 - \sum_{i=3}^{\Delta}(i - 1)N_i \left(\frac{2m}{n} - 2\right)^2 + \sum_{i=3}^{\Delta} N_i \left(i - \frac{2m}{n}\right)^2$$

$$= \sum_{i=3}^{\Delta} N_i(i - 1)(i - 2) - \frac{2(2m - n)(m - n)}{n}$$

and this gives $Var(G) = \frac{1}{n}\sum_{i=3}^{\Delta} N_i(i - 1)(i - 2) - \frac{2(2m-n)(m-n)}{n^2}$. ▫

**Corollary 28.** Let $G$ be an n-vertex connected unicyclic graph. Then $S(G) = 2\sum_{i=3}^{\Delta} N_i(i - 2)$, $Var(G) = \frac{1}{n}\sum_{i=3}^{\Delta} N_i(i - 1)(i - 2)$ and $n \times Var(G) - S(G) = \sum_{i=4}^{\Delta} N_i(i - 2)(i - 3)$.

**Proof.** Since $m = n$, the Theorem 27 gives the results. ▫

**Corollary 29.** Let $G(\neq C_n)$ be a connected unicyclic graph with $n$ vertices. Then $\Omega(G) = \frac{Var(G)}{S(G)} \geq \frac{1}{n}$, the equality holds if and only if $Ds(G) \subseteq \{1,2,3\}$.

**Proof.** By Corollary 28, $n \times Var(G) - S(G) \geq 0$, the equality holds if and only if $Ds(G) \subseteq \{1,2,3\}$. Thus, $\Omega(G) = \frac{Var(G)}{S(G)} \geq \frac{1}{n}$ and equality holds if and only if $Ds(G) \subseteq \{1,2,3\}$. ▫

**Corollary 30.** Let $G$ be a connected irregular graph with $n$ vertices and cyclomatic number $c$ such that $1 \leq c \leq \frac{n+2}{2}$. Then $\Omega(G) \geq \frac{1}{n} - \frac{2}{S(G)}$.

**Proof.** By Theorem 27, we have,

$$S(G) - n \times Var(G) = \frac{2(2m - n)(m - n)}{n} + \frac{2}{n}\sum_{i=3}^{\Delta} N_i(in - 2n) - \sum_{i=3}^{\Delta} N_i(i - 1)(i - 2)$$

$$\leq \frac{2(2m - n)(m - n)}{n} - \sum_{i=4}^{\Delta} N_i(i - 2)(i - 3) \leq \frac{(2m - n)(2m - 2n)}{n}$$

$$\leq \frac{(3n - n)(3n - 2n)}{n} = 2n.$$

Therefore, $\Omega(G) = \frac{Var(G)}{S(G)} \geq \frac{1}{n} - \frac{2}{S(G)}$. ▫

**Proposition 31.** There exist infinitely many tridegreed unicyclic graphs $G$ with degree set $Ds(G) = \{1, 2, \Delta \geq 3\}$ for which $S(G) = IRD(G) = \frac{2N_\Delta N_1}{N_\Delta + N_1}(\Delta - 1) = 2N_1$.

**Proof.** Consider an arbitrary bidegreed unicyclic graph $H$ with degree set $Ds(H) = \{\delta = 1, \Delta \geq 3\}$ and vertex number $n(H) = N_1 + N_\Delta$. Denote by $G$ the tridegreed unicyclic graph generated from $H$ by inserting $N_2$ vertices of degree 2 into $H$. The number of vertices of $G$ will be equal to $n(G) = N_1 + N_2 + N_\Delta$. By Corollary 28 the unicyclic graphs G and H have identical degree deviations, $S(G) = S(H) = 2N_\Delta(\Delta - 2)$. Additionally, by Theorem 21, $N_1 = N_\Delta(\Delta - 2)$. Therefore, $IRD(G) = IRD(H) = \frac{2N_\Delta N_1}{N_\Delta + N_1}(\Delta - 1) = \frac{2N_\Delta^2}{N_\Delta(\Delta-1)}(\Delta - 2)(\Delta - 1) = 2N_\Delta(\Delta - 2)$. Thus $S(G) = S(H) = IRD(G) = IRD(H) = 2N_1$. ☐

**Proposition 32.** Let $T$ be an $n \geq 2$ vertex tree which contains $N_1 \geq 2$ pendent vertices (vertices of degree 1). Then $S(T) = \frac{2(n-2)}{n}N_1$.

**Proof.** By Theorem 22, $S(T) = \frac{4(n-2)}{n} + \frac{2(n-2)}{n}\sum_{i=3}^{\Delta} N_i(i-2)$. Therefore, by Theorem 21, $S(T) = \frac{4(n-2)}{n} + \frac{2(n-2)}{n}(N_1 - 2) = \frac{2(n-2)}{n}N_1$. ☐

**Proposition 33.** Let $G$ be an $n \geq 3$ vertex unicyclic graph which contains $N_1 \geq 0$ pendent vertices. Then $S(G) = 2N_1$.

**Proof.** Since $m = n$, by Theorem 21 and Corollary 28, $S(G) = 2N_1 = 2\sum_{i=3}^{\Delta} N_i(i-2)$. ☐

Let $T_n$ and $U_n$ be a tree and an unicyclic graph on $n \geq 4$ vertices, respectively. If $N_1(T_n) = N_1(U_n)$, then by Propositions 32 and 33, $S(T_n) \to S(U_n)$ when $n \to \infty$.

**Proposition 34.** Let $G$ be a connected graph with $n$ vertices, $N_1 \geq 0$ pendent vertices and cyclomatic number $c$ such that $1 \leq c \leq \frac{n+2}{2}$. Then, $S(G) \leq 2(N_1 + 2c - 2)$, with equality if and only if $G$ is an unicyclic graph.

**Proof.** Since $m \geq n$, by Corollary 28 we have $S(G) = \frac{2}{n}\sum_{i=3}^{\Delta} N_i(in - 2m) \leq 2\sum_{i=3}^{\Delta}(i-2)N_i$ and equality holds if and only if $G$ is an unicyclic graph. Therefore, by Theorem 21 we have $S(G) \leq 2\sum_{i=3}^{\Delta}(i-2)N_i = 2(N_1 + 2c - 2)$ and equality holds if and only if $G$ is an unicyclic graph. ☐

**Proposition 35.** There are infinitely many $n$-vertex connected bidegreed graphs $G(\Delta, \Delta - 1)$. For such bidegreed graphs, $\Omega(G(\Delta, \Delta - 1)) = \frac{Var(G(\Delta,\Delta-1))}{S(G(\Delta,\Delta-1))} = \frac{1}{2n}$.

**Proof.** For $n \geq 3$, the path $P_n$ is a bidegreed graph such that $\Delta - \delta = 2 - 1 = 1$. Therefore, there are infinitely many $n$-vertex connected bidegreed graphs $G(\Delta, \Delta - 1)$. Suppose $G(\Delta, \Delta - 1)$ is a bidegreed graph. By Proposition 14, $\Omega(G(\Delta, \Delta - 1)) = \frac{Var(G(\Delta,\Delta-1))}{S(G(\Delta,\Delta-1))} = \frac{\Delta-(\Delta-1)}{2n} = \frac{1}{2n}$. ☐

## 8. Two Conjectures

Let $G$ be an arbitrary $n$-vertex and $m$-edge simple graph (connected or disconnected) with maximum degree $\Delta$ and minimum degree $\delta$. Our calculations on trees, unicyclic graphs and graphs with a small numbers of vertices suggest the following two conjectures:

**Conjecture 1.** It is conjectured that

$$S(G) \geq IRD(G) \quad \text{and} \quad Var(G) \geq \frac{1}{n^2} IRR(G) IRD(G).$$

The equalities hold if and only if $Ds(G) \subseteq \{\delta, \frac{2m}{n}, \Delta\}$.

**Conjecture 2.** It is conjectured that $\Omega(G) = \frac{Var(G)}{S(G)} \geq \frac{1}{2n}$.

By Proposition 35, there are infinitely many connected graphs such that $\Omega(G) = \frac{Var(G)}{S(G)} = \frac{1}{2n}$. It is important to note that sharp inequalities represented by above conjectures are valid for several tridegreed graphs.

**Example 5** As an example, consider the complete tripartite graph K$_{2,3,5}$ depicted in Fig.4.

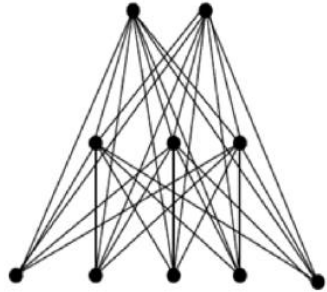

**Figure 4.** Complete tripartite graph K$_{2,3,5}$

Starting with the 10-vertex tridegreed graph K$_{2,3,5}$ characterized by degree sequence $(8^2, 7^3, 5^5)$ one obtains that $S(K_{2,3,5}) = 12 > 60/7 = IRD(K_{2,3,5})$. Because IRR(K$_{2,3,5}$)=15 and Var(K$_{2,3,5}$)=39/25, it follows that

$$1.56 = \frac{39}{25} = Var(K_{2,3,5}) = \frac{1}{n^2} IRR(K_{2,3,5}) IRD(K_{2,3,5}) = \frac{15}{100} \times \frac{60}{7} = \frac{9}{7} = 1.2857.$$

Moroever,

$$\Omega(K_{2,3,5}) = \frac{Var(K_{2,3,5})}{S(K_{2,3,5})} = \frac{39}{25 \times 12} = \frac{13}{25 \times 4} = \frac{13}{100} > \frac{5}{100} = \frac{1}{20} = \frac{1}{2n}.$$

## 9. Additional Considerations

Determination of extremal graphs with respect to irregularity measures $S(G)$ and $Var(G)$ belongs to the fundamental problems. In the research focused on the identification of extremal graphs the complete split graphs $CS(n,k)$ play a central role [11, 13, 14].

For complete split graphs one obtains that

$$S(CS(n,k)) = \frac{2k}{n}(n-k)(n-1-k) = \frac{2\delta}{n}(n-\delta)(n-1-\delta),$$

$$Var(CS(n,k)) = \frac{k}{n^2}(n-k)(n-1-k)^2 = \frac{\delta}{n^2}(n-\delta)(n-1-\delta)^2$$

Consequently, $\frac{1}{2n} \leq \Omega(CS(n,k)) = \frac{Var(CS(n,k))}{S(CS(n,k))} = \frac{n-k-1}{2n} = \frac{n-\delta-1}{2n} \leq \frac{1}{2}$.

Recently, it has been proved [14] that among all $n$-vertex connected graphs the maximal degree deviation is obtained for a particular subset of complete split graphs CS(n,k) where

- $k = n/3$, if $n$ is divisible by 3,
- $k = (n-1)/3$, if $n-1$ is divisible by 3,
- $k = (n-2)/3$ if $n-2$ is divisible by 3 and
- $k = (n+1)/3$ if $n+1$ is divisible by 3.

**An open problem:** Let $\Gamma(n,m)$ be the set of connected graphs with $n$ vertices and $m$ edges. Let $G$ be a connected irregular graph with $n$ vertices and m edges for which the degree deviation $S(G)$ is maximal with respect to any connected graph contained in set $\Gamma(n,m)$.

*Question:* Is it true that if an $n$-vertex and m edge connected graph $G$ included in $\Gamma(n,m)$ has a maximal degree deviation (i.e. S(G) is maximal with respect to all graphs contained in $\Gamma(n,m)$) then this graph G also has a maximal degree variance $Var(G)$?

Based on computational results it can be proved that the answer to the question is positive in some particular cases, but for the majority of cases the answer is negative. This observation is demonstrated on the following examples.

*Case A*: Consider the set $\Gamma(7,11)$ of connected irregular graphs with 7 vertices and 11 edges. In $\Gamma(7,11)$ there exist exactly 15 irregular graphs having universal vertices. Among them there are 14 graphs with one universal vertex, and only one 7-vertex graph with 2 universal vertices. This graph having the degree sequence $(6^2, 2^5)$ is isomorphic to the complete split graph $CS(7,2)$, and it is characterized by the following parameters: $N_6=2$, $N_2=5$, m=11, $M_1=92$ and $Var(CS(7,2)) = 160/49 = 3.2653$. It is easy to see that graph $CS(7,2)$ is the unique graph having a maximal degree variance among 15 irregular graphs belonging to set $\Gamma(7,11)$. Additionally, from previous considerations it follows that graph $CS(7,2)$ is the unique maximal graph in $\Gamma(7,11)$ with respect to the degree deviation for which $S(CS(7,2)) = 80/7$ is fulfilled. It is easy to check that $S(CS(7,2)) = \frac{2n}{\Delta-\delta} Var(CS(7,2)) = \frac{2\times 7}{6-2}\left(\frac{160}{49}\right) = \frac{80}{7}$ holds.

As we can conclude, among irregular graph with 7 vertices and 11 edges, the complete split graph $CS(7,2)$ is the unique extremal graph having maximal degree deviation and maximal degree variance, simultaneously.

*Case B*: Among 12-vertex connected graphs, the complete split graphs CS(12,4) with edge-number 38 has the maximal degree deviation $S(CS(12,4)) = 37.333$. Its degree variance is $Var(CS(12,4)) = 10.889$. Performing computer search it can be verifed that in the set $\Gamma(12,30)$ there exists a connected graph denoted by $H_{12,30}$ for which $Var(H_{12,30}) = 12$ and S($H_{12,30}$)=36 hold. Consequently, one obtains that $S(H_{12,30}) = 36 < 37.333 = S(CS(12,4))$ and $Var(H_{12,30}) = 12 > 10.889 = Var(CS(12,4)$. By simple computation it is easy to show that graph $H_{12,30}$ is isomorphic to the 12-vertex and 30-edge complete split graph $CS(12,3)$.